\documentclass[12pt,a4paper]{article}

\usepackage{mathptmx}
\usepackage{latexsym}
\usepackage{amssymb}
\usepackage{amscd}
\usepackage{amsmath}
\begin{document}
\title{Betti numbers of transversal monomial ideals}
\author{
Rahim Zaare-Nahandi \\[.5cm]
{\footnotesize Department of Mathematics, Statistics and Computer
Science,}\\ {\footnotesize University of
Tehran,}\\ {\footnotesize P.O. Box 14155-6455, Tehran, Iran}\\
{\footnotesize E-mail: rahimzn@khayam.ut.ac.ir}}
\date{}
\maketitle \thispagestyle{empty}
\begin{abstract}
In this paper, by a modification of a previously constructed minimal
free resolution for a transversal monomial ideal, the Betti numbers
of this ideal is explicitly computed. For convenient characteristics
of the ground field, up to a change of coordinates, the ideal of
$t$-minors of a generic pluri-circulant matrix is a transversal
monomial ideal . Using a Gr\"obner basis for this ideal, it is shown
that the initial ideal of a generic pluri-circulant matrix is a
stable monomial ideal when the matrix has two square blocks. By
means of the Eliahou-Kervair resolution, the Betti numbers of this
initial ideal is computed and it is proved that, for some
significant values of $t$, this ideal has the same Betti numbers as
the corresponding transversal monomial ideal. The ideals treated in
this paper, naturally arise in the
study of generic singularities of algebraic varieties. 
\end{abstract}
\section*{1. Introduction}
Let $S=k[y_{i,j(i)}: 1\leq i \leq n, 1\leq j(i)\leq b_i]$ be the
polynomial ring in $m=b_1+\cdots+b_n$ indeterminates over a field
$k$. Let ${\bf D}={\bf D}(b_1,b_2,\ldots,b_n)$ be the matrix
$$\left[
\begin{array}{cccccccccccc}
y_{11}& \cdots& y_{1b_1}&0&\cdots&0&0&\cdots&0&0&\cdots&0\\
0&\cdots&0&y_{21} &\cdots&y_{2b_2}&0&\cdots&0&0&\cdots&0\\
\vdots&\vdots&\vdots&
\vdots&\vdots&\vdots&\vdots&\vdots&\vdots&\vdots&\vdots&\vdots\\
0&\cdots&0&0&\cdots&0&0&\cdots&0&y_{n1}&\cdots&y_{nb_n}
\end{array}\right].$$
Let $I_t({\bf D})\in S$ be the ideal generated by $t$-minors of
${\bf D}$. This is a square-free monomial ideal and is called a {\it
transversal monomial ideal}. In (\cite{ZZ}, \S 3) the Hilbert series
of $S/I_t({\bf D})$ has been computed by means of the simplicial
complex associated to $I_t({\bf D})$ and its minimal free resolution
has been constructed. In this paper we outline a modification of
this resolution and compute the Betti numbers of
$I_t({\bf D})$.\\

Let $R=k[x_{ij}: 1\leq i\leq n, 1\leq j\leq b]$ be the polynomial
ring in $nb$ indeterminates over $k$ and let ${\bf
P}=\left[\begin{matrix}M_1& M_2&\cdots &M_b\end{matrix}\right]$ be a
{\it generic pluri-circulant matrix} where $M_j$ is the generic
circulant matrix with the first row
$\left(\begin{matrix}x_{1j}&x_{2j}&.&.&.&x_{nj} \end{matrix}
\right).$ The ideal generated by $t$-minors of ${\bf P}$ has been
considered in \cite{Z}. Under some hypothesis on $k$, this ideal is
closely related to a certain transversal monomial ideal. In fact, if
$k$ possesses the $n$th roots of unity and ${\rm char}(k)\nmid n$
then over such ground field, up to a linear change of coordinates,
the matrix ${\bf P}$ is equivalent to the matrix ${\bf D}$ with
$b_1=b_2=\cdots=b_n=b$ (\cite{ZZ}, \S 4). On the other hand, with
the same assumptions on the ground field, it has been proved that
for a suitable monomial order on $R$, certain set of $t$-minors of
the first $t$ rows of ${\bf P}$ forms a Gr\"obner basis for
$I_t({\bf P})$ and its initial ideal $J_t$ has been computed (see
\cite{SZ}, \S 3). For a filed of arbitrary characteristic, such a
result is known only for $t=n, n-1$ (\cite{Z}, \S 5). However, the
monomial ideal $J_t$ can be studied in its own.  We show that for
$b=2$, $J_t$ is a {\it stable monomial ideal}. This class of
monomial ideals have been introduced and studied by Eliahou and
Kervaire \cite{EK}. Using the Eliahou-Kervaire resolution for stable
monomial ideals, we compute the Betti numbers of $J_t$. For $t=n$,
$n-1$ and $n-2$, we prove that the Betti numbers of $J_t$ are equal
to the corresponding Betti numbers of $I_t({\bf D})$. These
equalities are not immediate and require some unexpected
combinatorics.\\

The ideals treated here, naturally arise in the study of the local
equations of generic singularities of algebraic varieties (see
\cite{SZ}, \cite{Z'}).
\section*{2. The minimal free resolution and Betti numbers of $I_t({\bf
D})$} The notation employed for description of the minimal free
resolution of $I_t({\bf D})$ in (\cite{ZZ} \S 3) can be modified to
make this resolution more accessible. This is the first task of the
section. Using this setting, we compute the Betti numbers of
$I_t({\bf D})$ explicitly. In the special case $m=n$, the
modification allows one to define a structure of a graded
differential algebra on the resolution. \\

Let $V$ be a free $S$-module of rank $m=b_1+\cdots+b_n$ generated by
symbols $e_{i,j(i)}$ in one-to-one correspondence with the
indeterminates $y_{i,j(i)}$. Let $W$ be the free $S$-module of rank
$n$ generated by symbols $\epsilon_1, \cdots, \epsilon_n$ in
one-to-one correspondence with rows of ${\bf D}$. For $q=0,
1,\cdots, m-t-1$, let $E_q= \bigoplus_{p=1}^{q+1} \wedge^pW$ and let
$C_q \subset (\wedge^{t+q}V)\otimes E_q$ be the free $S$-module
generated by the basis elements
$$
e_{i_1,j_1(i_1)}\wedge \cdots \wedge
e_{i_1,j_{r_1}(i_1)}\wedge\cdots \wedge e_{i_s,j_1(i_s)}\wedge
\cdots \wedge e_{i_s,j_{r_s}(i_s)}\otimes \epsilon_{i_{k_1}}\wedge
\cdots \wedge \epsilon_{i_{k_{s-t+1}}}
$$
where
$$t\leq s\leq n, \ \ \ r_1,\cdots, r_s \geq 1, \ \ {\rm
and}, \ \ \ 1\leq k_1<\cdots < k_{s-t+1}\leq s.$$ As for the
elements of the wedge product, we adopt the usual convention on the
order of vectors $e_{i,j(i)}$ to appear in the lexicographic order
of their indices, i.e.,
$$1\leq i_1< \cdots <i_s \leq n,$$
and
$$1\leq j_u(i_v) < j_{u+1}(i_v) \leq
b_{i_v}, \forall u, v.$$

\noindent Clearly, $r_1+\cdots +r_s =t+q,$ and $t\leq s\leq$
Min$\{n, t+q\}$.\\

For simplicity, we may drop the first subscripts in each
$e_{i,j(i)}$ and denote the basis elements of $C_q$ by
$$
e_{j_1(i_1)}\wedge \cdots \wedge e_{j_{r_1}(i_1)}\wedge\cdots \wedge
e_{j_1(i_s)}\wedge \cdots \wedge e_{j_{r_s}(i_s)}\otimes
\epsilon_{i_{k_1}}\wedge \cdots \wedge \epsilon_{i_{k_{s-t+1}}}.
$$

To keep a reference of the above basis elements, it may be helpful
to replace $y_{i,j(i)}$ with $e_{j(i)}$ in the matrix ${\bf D}$.
Then the above basis elements are obtained by selecting the rows
$i_1, \cdots, i_s$ in the resulting matrix and choosing $r_1$
nonzero entries on the $i_1$st row, ... , and, $r_s$ nonzero entries
on the $i_s$th row, and finally, specifying further rows
$i_{k_1}, \cdots, i_{k_{s-t+1}}$ among the selected rows.\\

We fix some notation to use in the sequel. Let $$\xi =
e_{j_1(i_1)}\wedge \cdots \wedge e_{j_{r_1}(i_1)}\wedge\cdots \wedge
e_{j_1(i_s)}\wedge \cdots \wedge e_{j_{r_s}(i_s)})\in
\wedge^{t+q}V,$$
$$\delta = \epsilon_{i_{k_1}}\wedge \cdots \wedge
\epsilon_{i_{k_{s-t+1}}}\in E_q .$$ Furthermore, for
$v\in\{1,\cdots, s\}$, and $h\in\{1,\cdots,r_v\}$, let
$$\xi_{\widehat{e_{j_h(i_v)}}}= e_{j_1(i_1)}\wedge \cdots\wedge
e_{j_1(i_v)}\wedge \cdots \wedge \widehat{e_{j_h(i_v)}}\wedge \cdots
\wedge e_{j_{r_v}(i_v)}\wedge \cdots \wedge e_{j_{r_s}(i_s)},$$
where the hat sign over a term means that this term is to be
omitted. For $w\in\{1,\cdots,s-t+1\}$ with $r_{k_w}=1,$ let
$$\xi_{\widehat{e_{j_1(i_{k_w})}}}\otimes\delta_{\widehat{\epsilon_{i_
{k_w}}}} = e_{j_1(i_1)}\wedge \cdots \wedge
\widehat{e_{j_1(i_{k_w})}}\wedge\cdots \wedge
e_{j_{r_s}(i_s)}\otimes \epsilon_{i_{k_1}}\wedge
\cdots\wedge\widehat{\epsilon_{i_{k_w}}}\wedge\cdots \wedge
\epsilon_{i_{k_{s-t+1}}}.$$

Similar to (\cite{EN}, (3.2)), for $i_v\in \{i_1,\cdots,i_s\}$, let
$\Delta_v$ be the ``{\it differentiation by the $i_v$-th row of
${\bf D}$}''. I.e., for $v$ with $r_v\geq 2,$ let
$$\Delta_v(\xi)=
\sum_{h=1}^{r_v}(-1)^{r_1+\cdots+r_{v-1}+h+1} y_{i_v,j_h(i_v)}
\xi_{\widehat{e_{j_h(i_v)}}}.$$ For $r_v=1$, $\Delta_v(\xi)$ is
defined to be zero. Then as in (\cite{EN}, (3.3)),
$$\Delta_u\circ \Delta_v +\Delta_v\circ \Delta_u =0.\eqno(2.1)$$

Furthermore, for $1\leq k_1, \cdots, k_{s-t+1}\leq s$, define
$\Lambda = \Lambda_{k_1,\cdots,k_{s-t+1}}$ as the ``{\it
differentiation by diagonal entries}''. For $s-t+1 \geq 2$ let
$$\Lambda(\xi\otimes \delta)=
\sum_{w^*}
(-1)^{\ell(w)+1}y_{i_{k_w},j_1(i_{k_w})}\xi_{\widehat{e_{j_1(i_{k_w})}}}
\otimes\delta_{\widehat{\epsilon_{i_{k_w}}}}$$ where the asterisk
sign over $w$ means that we sum only over those values of $w$ for
which $r_{k_w}=1,$ and $\ell(w)=\#\{k_v| v\leq w, r_{k_v}=1\}.$ For
$s-t+1=1$, $\Lambda(\xi\otimes \delta)$ is defined to be zero. It
then follows that
$$\Lambda^2(\xi\otimes \delta)=0. \eqno(2.2)$$

For $q\geq 0$ the boundary map $d:C_{q+1} \longrightarrow C_q$ can
now be defined by
$$d(\xi \otimes \delta)= \sum
_{v^*}\Delta_v(\xi)\otimes \delta + \Lambda(\xi\otimes \delta),$$
where the asterisk sign over $v$ in the summation means that we sum
only over those values of $v$ for which $r_v\geq 2.$

By (2.1), (2.2) and some straightforward computation, it follows
that $C_{\mbox{\tiny$\bullet$}}$ is indeed a complex. However, in
general, the complex $C_{\mbox{\tiny$\bullet$}}$ does not lead to a
free resolution for $I_t({\bf D})$. We construct a subcomplex
$K_{\mbox{\tiny$\bullet$}}\subset C_{\mbox{\tiny$\bullet$}}$ such
that an augmentation of the quotient complex $
C_{\mbox{\tiny$\bullet$}}/K_{\mbox{\tiny$\bullet$}}$ is the minimal
free resolution for $I_t({\bf D}).$

For $q\geq 0,$ let $K_q\subset C_q$ be the submodule generated by
all expressions
$$\xi\otimes(
\sum_{p=1}^{s-t+2} (-1)^p \epsilon_{i_{k_1}}\wedge \cdots
\wedge\widehat{\epsilon_{i_{k_p}}}\wedge\cdots \wedge
\epsilon_{i_{k_{s-t+2}}})$$ for some $k_1, \cdots, k_{s-t+2}$ with
$1\leq k_1<\cdots <k_{s-t+1}<k_{s-t+2}=s,$ where as above, $$\xi =
e_{j_1(i_1)}\wedge \cdots \wedge e_{j_{r_1}(i_1)}\wedge\cdots \wedge
e_{j_1(i_s)}\wedge \cdots \wedge e_{j_{r_s}(i_s)}.$$ It can be
checked that $K_{\mbox{\tiny$\bullet$}}$ is a sub-complex of
$C_{\mbox{\tiny$\bullet$}},$ i.e., $d(K_{q+1})\subset d(K_q).$ Let
$L_q=C_q/K_q$ for $q=0, \cdots, m-t-1$. While among the summands of
any of the above expressions there is only one $\epsilon$ without an
index equal to $i_s$, we will consider the representatives
$$e_{j_1(i_1)}\wedge \cdots \wedge e_{j_{r_1}(i_1)}\wedge\cdots
\wedge e_{j_1(i_s)}\wedge \cdots \wedge e_{j_{r_s}(i_s)})\otimes
\epsilon_{i_{k_1}}\wedge \cdots \wedge \epsilon_{i_{k_{s-t}}}\wedge
\epsilon_{i_s}, \eqno(2.3)$$ as basis elements of $L_q$. In
particular $L_q$ is a free $S$-module. Although one may ignore
writing the index $i_s,$ it may be kept to signify the action of the
differentiations by diagonal entries.

Finally, the augmentation map $d:L_0 \longrightarrow I_t({\bf D})$
is defined as the {\it determinant map}, i.e.,

$$d(e_{j_1(i_1)}\wedge \cdots \wedge e_{j_1(i_t)}\otimes \epsilon_{i_t})
=y_{i_1j_1(i_1)}\cdots y_{i_tj_1(i_t)}.$$\\

The main result of \cite{ZZ} may now be stated. For the proof we
refer to (\cite{ZZ}, Theorem 3.1).
\paragraph{2.1. Theorem.} The complex $L_{\mbox{\tiny$\bullet$}}$
with the induced boundary maps is the minimal free resolution for
$I_t({\bf D})$.
\paragraph{2.2. Remark.} The complex $L_{\mbox{\tiny$\bullet$}}$
is clearly linear. However, the linearity also follows via
properties of the simplicial complex associated to $I_t({\bf D})$
(see \cite{ZZ}, Proposition 2.1). In fact, $I_t({\bf D})$ is {\it
weakly polymatroidal}, and hence, it has linear quotients and in
particular, it has linear resolution (see \cite{KH}, Theorem 1.4).\\

The natural multiplication
$$(\xi \otimes \delta).(\zeta \otimes \tau) = (\xi \wedge
\zeta)\otimes (\delta\wedge \tau), \eqno (2.4)$$ is not well-defined
on the complex $C_{\mbox{\tiny$\bullet$}}$ to turn it into a
differential algebra unless $t=n$, or, $b_1=\cdots =b_n =1$. Even
for $t=n$, the Leibnitz formula fails. However, in the latter case,
$C_{\mbox{\tiny$\bullet$}}$ is a graded differential algebra under
the above multiplication. Thus we may state the following which
should have been well-known.
\paragraph{2.3. Corollary.}
For $m=n$ the complex $C_{\mbox{\tiny$\bullet$}}$ is a graded
differential algebra under the multiplication (2.4) and the complex
$K_{\mbox{\tiny$\bullet$}}$ is a homogenous ideal in
$C_{\mbox{\tiny$\bullet$}}$ and hence $L_{\mbox{\tiny$\bullet$}}$
inherits the structure of a graded differential algebra as the
minimal free resolution of the ideal generated by all square-free
monomials of degree $t$ in $n$ indeterminates. In this case, with
the exception of the augmentation map, the
boundary maps descend to $\Lambda$.\\

We now return to the general case. Although when the minimal free
resolution is known, the Betti numbers are theoretically available,
explicit formulation of these numbers provide finer information. We
now pursue on the computation of the Betti numbers of $I_t({\bf
D})$.
\paragraph{2.4. Proposition.} With the notations as the above,
$$\beta_q(I_t({\bf D})) = \sum_{s=t}^{{\rm Min}\{t+q,n\}} {s-1\choose t-
1} \ \sum_{1\leq i_1<\cdots < i_s\leq n} \ \ \ \sum_{\begin{matrix}
r_1+\cdots+r_s=t+q \cr r_1,\cdots,r_s\geq 1\end{matrix}}
{b_{i_1}\choose r_1}\cdots{b_{i_s}\choose r_s}.$$ For the case
$b_1=\cdots=b_n=b$,
$$\beta_q(I_t({\bf D})) = \sum_{s=t}^{t+q}{s-1\choose t-1}{n\choose s}
\ \ \ \sum_{\begin{matrix} r_1+\cdots+r_s=t+q \cr r_1,\cdots,r_s\geq
1\end{matrix}} {b\choose r_1}\cdots{b\choose r_s}.$$ For
$b_1=\cdots=b_n=2$,
$$\beta_q(I_t({\bf D})) = \sum_{s=t}
^{t+q} {s-1\choose t-1}{n\choose s}{s\choose t+q-s}2^{2s-t-q}.$$ The
usual conventions ${\alpha\choose \beta}=0$ for $\alpha < \beta$,
and ${\alpha\choose 0} = 1$ for $\alpha \geq 0$, are to be adopted
in these formulas. Thus, for example, the precise lower bound and
upper bound of the last summation would be $s={\rm Max}\{t,
\lceil\frac{t+q}{2}\rceil\}$ and $s={\rm Min}\{t+q,n\}$,
respectively.
\paragraph{Proof.} We use the expressions (2.3) for the basis of $L_q$.
In the first formula, ${s-1\choose t-1}$ is the number of choices
for $k_1,\cdots, k_{s-t}$, the second summation is for the choices
of $i_1,\cdots,i_s$ and the last summation is for the number of
choices of $r_1, \cdots, r_s.$ The second equality is immediate. For
the last equality, while for each $v\in \{1,\cdots,s\}$, $r_v=1 \
{\rm or} \ 2$, the binomial coefficient counts the number of cases
and the power of $2$ is the number of $r_v$'s with
$r_v=1$.$\hfill\square$
\paragraph{2.5. Remark.}
It is important to emphasize the condition $r_1,\cdots,r_s\geq 1$ in
the first and the second formulas in the above proposition.
Otherwise, the last summation would simplify using the following
so-called {\it generalized Vandermonde convolution} \cite{GKP}:
$$\sum_{r_1+\cdots+r_s=t+q} {b_{i_1}\choose r_1}\cdots
{b_{i_s}\choose r_s}= {b_{i_1}+\cdots + b_{i_s}\choose
r_1+\cdots+r_s}.$$
\section*{1. Betti numbers of the initial ideal of the ideal of
$t$-minors of generic pluri-circulant matrices}

Generic pluri-circulant matrices and their ideals of $t$-minors
arise naturally in the study of generic projections in algebraic
geometry \cite{SZ}. These ideals are closely related to transversal
monomial ideals. In this section we recall some results on the
ideals of $t$-minors of generic pluri-circulant matrices and then we
determine Betti numbers of the initial ideals of these ideals of
minors and compare
them with the Betti numbers computed in the previous section.\\

Let $R=k[x_{ij}: 1\leq i\leq n, 1\leq j\leq b]$ be the polynomial
ring in $nb$ indeterminates over $k$ and let ${\bf
P}=\left[\begin{matrix}M_1& M_2&\cdots &M_b\end{matrix}\right]$ be a
generic pluri-circulant matrix where
$$M_j= \left[\begin{matrix}x_{1j}&x_{2j}&.&.&.&x_{nj} \cr
x_{nj}&x_{1j}&.&.&.&x_{n-1,j} \cr \vdots
&\vdots&\vdots&\vdots&\vdots&\vdots\cr
x_{2j}&x_{3j}&.&.&.&x_{1j}\end{matrix} \right]$$ is a generic
circulant matrix. Let $k$ possess the $n$th roots of unity and ${\rm
char}(k)\nmid n$. Then, under a linear change of variables in $R$,
$I_t({\bf P})$ converts to the ideal  $I_t({\bf D})$ in $S$
considered in the previous section with $b_1=b_2=\cdots=b_n=b$ (see
\cite{ZZ}, $\S$ 4). Let
$${\bf T} = \left [\begin{matrix}
T_1& T_2& \cdots & T_b\end{matrix}\right ],$$ where
$$T_j=\left [\begin{matrix}x_{1j}&x_{2j}&\cdots&.&.&\cdots&x_{nj}\cr
0&x_{1j}&\cdots&.&.&\cdots&x_{n-1,j}\cr\vdots&\vdots&\vdots&\vdots
&\vdots&\vdots&\vdots
\cr0&0&\cdots&0&x_{1j}&\cdots&x_{n-t+1,j}\end{matrix}\right], \ \
j=1,2,\cdots,b.$$ Let $J_t$ be the ideal in $R$ generated by
products of the entries of the main diagonals of ${\bf T}$. The
ideal $J_t$ is $Q$-primary where $Q$ is the prime ideal generated by
the indeterminates on the last row of ${\bf T}$. It is known that,
under the above assumptions on the ground field, for a suitable
monomial order on $R$, the set of $t$-minors of the first $t$ rows
of ${\bf P}$ whose main diagonals correspond to the main non-zero
diagonals of ${\bf T}$, forms a Gr\"obner basis for $I_t({\bf P})$
and $J_t$ is its initial ideal (\cite{SZ}, Theorem 3.3). For
arbitrary filed $k$, this result is only known for $b=2, t=n, n-1$
(\cite{Z}, Theorem 5.4).  We show that for $b=2$, $J_t$ is a {\it
stable monomial ideal}. Using the Eliahou-Kervaire resolution for
stable monomial ideals \cite{EK}, we prove that all Betti numbers of
$J_t$ and $I_t({\bf D})$ are equal at least for $t=n, n-1, n-2$. In
general, Betti numbers of the initial ideal of a given ideal, only
give upper bounds
for Betti numbers of the the original ideal.\\

Recall that a monomial ideal $I\subset k[z_1,\cdots,z_n]$ is said to
be {\it stable} if for every monomial $w\in I$ and index $i < m =
{\rm max}(w)$, the monomial $z_iw/z_m$ again belongs to $I$, where
${\rm max}(w)$ denotes the largest index of the variables dividing
$w$. Let $G(I)$ be the unique minimal generating set of $I$
consisting of monomials. Note that $I$ is stable if and only if the
above condition holds for every $w\in G(I)$. Clearly, no nontrivial
square-free monomial ideal is stable. In particular, $I_t({\bf D})$
is not a stable monomial ideal. For $b=2$, to simplify the notation
we use $d=n-t+1$, and we consider the re-indexing of indeterminates
in the ring $R$ such that the matrices $T_1$ and $T_2$ turn to the
following forms, respectively:
$$T'_1= \left [\begin{matrix}z_1&z_2&\cdots&z_d&z_{2d+1}&z_{2d+2}&\cdots
&z_{2d+t-1}\cr
0&z_1&\cdots&z_{d-1}&z_d&z_{2d+1}&\cdots&z_{2d+t-2}\cr
\vdots&\vdots&\vdots&\vdots&\vdots&\vdots&\vdots &\vdots
\cr0&\cdots&0&z_1&z_2&\cdots&z_{d-1}&z_d\end{matrix}\right],$$
$$T'_2=\left[\begin{matrix} z_{d+1}&z_{d+2}&\cdots&z_{2d}&
*&* &\cdots&*\cr 0& z_{d+1}&\cdots&z_{2d-1}&z_{2d}&*&\cdots &*\cr
\vdots&\vdots&\vdots&\vdots&\vdots&\vdots&\vdots &\vdots\cr
0&\cdots&0&z_{d+1}&z_{d+2}&\cdots&z_{2d-1}&z_{2d}
\end{matrix}\right].$$
In other words, we first re-index the indeterminates on the last row
of $\left[\begin{matrix} T_1&T_2\end{matrix}\right]$ so that the
nonzero entries on the last row are linearly ordered. Then the
remaining indeterminates on the first row of this matrix are
re-indexed linearly. Since the last $t-1$ indeterminates on the
first row do not appear in $G(J_t)$, they are replaced by $*$'s.
With this new indexing, the ideal $J_t$ is stable.
\paragraph{3.1. Lemma.} For $b=2$, with the above re-indexing
of the indeterminates, $J_t$ is a stable monomial ideal. In
particular, the associated Eliahou-Kervaire complex provides a
linear minimal free resolution for $J_t$ equipped with a structure
of graded differential algebra.
\paragraph{Proof.} Recall that $G(J_t)$ consists of the products of
entries of the main diagonals of $\left[\begin{matrix}
T'_1&T'_2\end{matrix}\right]$ . Observe that every monomial in
$G(J_t)$ has a unique representation in the form
$$z_{i_1}\cdots z_{i_q}z_{j_1}\cdots z_{j_r} z_{k_1}\cdots
z_{k_s}$$
with\\
$$1\leq i_1 \leq \cdots \leq i_q\leq d,$$
$$2d+1\leq j_1 \leq \cdots \leq j_r\leq 2d+t-1-q,$$
$$d+1\leq k_1 \leq \cdots \leq k_s\leq 2d,$$
$$q+r+s=t.$$
Conversely, any such representation identifies a unique monomial in
$G(J_t)$. Let $w\in G(J_t)$ and let $w=z_{i_1}\cdots
z_{i_q}z_{j_1}\cdots z_{j_r} z_{k_1}\cdots z_{k_s}$ be its unique
representation. Since all monomials of degree $t$ in $z_\ell$ with
$1\leq \ell \leq 2d$ belong to $J_t$, we need to check the stability
condition for the case $r\geq 1$. Then the maximum index of $w$ is
$j_r$. We need to show that $w'=\frac{z_\ell w}{z_{j_r}}\in G(J_t)$
for all $1\leq \ell < j_r$. This can be checked directly. In fact,
if $i_\tau \leq \ell \leq i_{\tau +1}$, then $w'=z_{i_1}\cdots
z_{i_\tau}z_\ell z_{i_{\tau +1}}\cdots z_{i_q}z_{j_1}\cdots
z_{j_{r-1}} z_{k_1}\cdots z_{k_s}\in G(J_t).$ If $j_\tau \leq \ell
\leq j_{\tau +1}$, then $w'=z_{i_1}\cdots z_{i_q}z_{j_1}\cdots
z_{j_\tau}z_\ell z_{j_{\tau +1}}\cdots z_{j_{r-1}} z_{k_1}\cdots
z_{k_s}\in G(J_t).$ If $k_\tau \leq \ell \leq k_{\tau +1}$, then
$w'=z_{i_1}\cdots z_{i_q}z_{j_1}\cdots z_{j_{r-1}} z_{k_1}\cdots
z_{k_\tau}z_\ell z_{k_{\tau +1}}\cdots z_{k_s}\in G(J_t).$ The last
claim follows from (\cite{EK}, \S 2 Theorem 2.1 and Remark 1).
$\hfill\square$
\paragraph{3.2. Remark.} For $b>2$, the ideal $J_t$ is not a
stable monomial ideal as it can be inspected for $b=3$, $n=t=3$. On
the other hand, even for $b=2$, the ideal $J_t$ is not Borel fixed
(see the definition in \cite{BS} or
\cite{EK}), as it can be checked for $n=t=5$.\\

Recall that by (\cite{EK}, \S 3), the Betti numbers of a stable
monomial idea $I$ is given by
$$\beta_q(I)=\sum_{w\in G(I)}{{\rm max}(w)- 1 \choose q},$$
where $G(I)$ is the minimal generating set of $I$. We will use this
result to compute the Betti numbers of $J_t$ explicitly.
\paragraph{3.3. Lemma.} For $b=2$, let $\nu_\ell$ be the number
of monomials in $G(J_t)$ with largest index $\ell$ and let
$d=n-t+1$.
Then\\
(a) For $1\leq \ell \leq 2d$, $$\nu_\ell = {t+ \ell -2\choose
t-1}.$$ (b) For $2d+1\leq \ell \leq 2d+t-1$,
$$\nu_\ell = \sum_{k=0}^{t-j-1} {n-t+j+k\choose k}{n-k-1 \choose
n-t}= \sum_{\tau=1}^{n-t+1} {n+\tau -1\choose t-j-1}{n-t+j-\tau
\choose j-1}$$
$$= \sum_{\tau=1}^{n-t+1} {n\choose t-j-\tau}{n-t+j
\choose j+\tau-1},$$
 where $j=\ell - 2d$.\\
(c) For $q=0,1,\cdots, 2n-t+1$,
$$\beta_q(J_t)=\sum_{\ell=q+1}^{2n-t+1}{\ell-1\choose
q}\nu_\ell,$$ where $\nu_\ell$ is given in (a) and (b).
\paragraph{3.4. Remark.} Observe that the first equality for $\nu_\ell$
in (b) also makes sense for $d+1\leq \ell \leq 2d$ and it reduces to
(a). This follows from the well-known identity
$$\sum_{k=a-m}^{c-n} {a+k\choose m}{c-k \choose n}= {a+c+1\choose
m+n+1}$$ for all non-negative integers $a, m, c, n$ with $a\geq m$
and $c\geq n$. However, the second and third equality in (b) is
valid only for $2d+1\leq \ell \leq 2d+t-1$.
\paragraph{Proof of Lemma 3.3.} The
ideal $J_t$ contains all monomials of degree $t$ in $z_1, \cdots,
z_{2d}$. Thus for $1\leq \ell \leq 2d$, a monomials in $G(J_t)$ with
largest index $\ell$ is of the form $wz_\ell$ where $w$ is any
monomial of degree $t-1$ in $z_1,\cdots, z_\ell$. This settles (a).
For $2d+1\leq \ell \leq 2d+t-1$ a monomial $w\in G(J_t)$ with
largest index $\ell=2d+j$ can be uniquely written in the form
$w=w_1z_\ell w_2$, where $w_1$  varies in the set of all monomials
of degree $k$ in $z_1,\cdots, z_d,z_{d+1}\cdots, z_\ell$ and $w_2$
ranges over all monomials of degree $t-k-1$ in
$z_{d+1},\cdots,z_{2d}$, for $k=0,\cdots, t-j-1$. This proves the
first equality in (b). To prove the second equality, we employ
another method to count the same monomials. Using the index
configuration in $\left[\begin{matrix}
T'_1&T'_2\end{matrix}\right]$, for fixed $j$, $1\leq j=\ell- 2d\leq
t-1$, any such monomial can be uniquely written in the form
$w=u_1z_\ell z_{d+\tau}u_2$, where $u_1$ varies over the set of all
monomials of degree $t-j-1$ in $z_1,\cdots, z_d,z_{2d+1}\cdots,
z_\ell, z_{d+1}, \cdots, z_{d+j+\tau}$ and $u_2$ ranges over all
monomials of degree $j-1$ in $z_{d+\tau},\cdots,z_{2d}$, for
$\tau=1,\cdots, d$. In fact, $u_1$ is a product any diagonal entries
of
$$\left [\begin{matrix}z_1&\cdots&z_d&z_{2d+1}&\cdots&
z_\ell&z_{d+1}&\cdots &z_{d+j+\tau}&\cdots\cr
\vdots&\vdots&\vdots&\vdots
&\vdots&\vdots&\vdots&\vdots&\vdots&\vdots\cr
\cdots&z_1&\cdots&z_d&z_{2d+1}&\cdots&z_\ell&z_{d+1}&\cdots
&z_{d+j+\tau}\end{matrix}\right]$$ as a submatrix of
$\left[\begin{matrix} T'_1&T'_2\end{matrix}\right]$ with $t-j-1$
rows, and $u_2$ is a product of any diagonal entries of
$$\left [\begin{matrix}z_\tau&\cdots&z_{2d}&\cdots\cr
\vdots&\vdots&\vdots&\vdots
 \cr
\cdots& z_\tau&\cdots&z_{2d}
\end{matrix}\right]$$
a submarix of $T'_2$ wit $j-1$ rows. This clarifies the second
equality in (b). The third equality in (b) follows from the second
equality as a combinatorial identity, or, by using the identity
$(3.2)$ below. The assertion (c) is just the formula for the Betti
numbers of a stable monomial ideal. $\hfill\square$
\paragraph{3.5. Remark.} While for $b=2$ the minimal free resolution of
$J_t$ has a natural structure of graded differential algebra, the
minimal free resolution of $I_t({\bf D})$ and hence that of
$I_t({\bf P})$ has no such natural structure as explained prior to
Corollary 2.3. In particular, the minimal free resolution of
$I_t({\bf P})$ is not a ``natural lifting'' of the Eliahou-Kervaire
resolution of $J_t$. More importantly, if ${\rm char}(k) | n$ and no
minimal free resolution for $I_t({\bf P})$ is known, the minimal
free resolution of $J_t$ does not naturally lift to the minimal free
resolution of $I_t({\bf P})$. This is contrary to what one could
have hoped, since, at least for $t=n, n-1$, $J_t= {\rm in}(I_t({\bf
P}))$. However, regardless of characteristic of the ground field,
these three ideals have the
same Betti numbers.\\

The following result should be true for all $1\leq t\leq n$. We only
prove it for $t=n, n-1, n-2$. Although the same procedure works for
any specific value of $t$, we are not able to provide a unified
proof for arbitrary $t$. We will replace $d$ with $n-t+1$.
\paragraph{3.6. Theorem.} For $b=2$, $t=n, n-1, n-2$, the ideals $J_t$,
$I_t({\bf D})$ and $I_t({\bf P})$ have equal Betti numbers, i.e.,
$$\beta_q(J_t) = \beta_q(I_t({\bf D})) = \beta_q(I_t({\bf P})),$$
for all $q=0,\cdots, 2n-t$.
\paragraph{Proof.} The last equality is clear due to the
explanations at the beginning of this section. Thus we prove the
first equality.\\
For $t=n$, the claim is rather straightforward. Indeed, by Lemma
3.3,
$$\beta_q(J_n)=\sum_{\ell=q+1}^{2}{\ell-1\choose
q}{n+ \ell -2\choose n-1} + \sum_{j=q-1}^{n-1}{j+1\choose
q}{n\choose j+1}.$$ For $q\geq 2$, the first sum is zero. By
Proposition 2.4,
$$\beta_q(I_t({\bf D}))= {n\choose q} 2^{n-q}.$$ Thus for $q\geq 2$
the equality $\beta_q(J_n)=\beta_q(I_t({\bf D}))$ is just a
well-known combinatorial identities. For $q=0,1$, the proof is
similar.\\
For $t=n-1, n-2$ to settle the equality $\beta_q(J_t)=
\beta_q(I_t({\bf D}))$ we try to write both sides as
$\mathbb{Z}[q]$-linear combinations of ${n\choose q -\alpha}
2^{n-q+\alpha}$ for $\alpha =0,\pm 1, \pm 2$. The main combinatorial
identities to be employed are $${i\choose \tau}{i-\tau\choose
q-\tau} = {q\choose \tau}{i\choose q}, \ \ \ \eqno (3.1)$$
$${n+\rho\choose q}= \sum_{i=0}^\rho {\rho\choose i}{n\choose
q-i} \ \ \   \eqno (3.2)$$ and
$$\sum_{i=q}^n{i\choose q}{n\choose i}= {n\choose q} 2^{n-q}. \ \ \ \eqno (3.3)$$
For $t=n-1,$ by (c), (a) and the last equality of (b) in Lemma 3.3),
$$\beta_q(J_{n-1})=\sum_{\ell=q+1}^{4}{\ell-1\choose q}{n+ \ell
-3\choose n-2} + \sum_{j=q-3}^{n-2}{j+3\choose q}\sum_{\tau =1}^2
{n\choose j+\tau +1}{j+1\choose 2-\tau}.$$ We treat the case $q \geq
4$ so that the first sum is zero. For $0\leq q\leq 3$, similar
computation works where the first sum recovers the missing quantity
expected for the required equality. For $q \geq 4,$ using $i=j+3$ we
get
$$\beta_q(J_{n-1})= \sum_{i=q}^{n+1}(i-2){i\choose q}{n\choose i-1} +
\sum_{i=q}^n{i\choose q}{n\choose i}.$$ By $(3.1)$ the first sum in
$\beta_q(J_{n-1})$ reduces to
$$\sum_{i=q}^{n+1}(i+1){i\choose q}{n\choose
i-1}-3\sum_{i=q}^{n+1}{i\choose q}{n\choose i-1}$$ $$=
\sum_{i=q}^{n+1}(q+1){i+1\choose q+1}{n\choose i-1}-
3\sum_{i=q+1}^{n+1}{i-1\choose q}{n\choose i-1}
-3\sum_{i=q}^{n+1}{i-1\choose q-1}{n\choose i-1}.$$ Using $(3.2)$
and $(3.3)$ we have $$\sum_{i=q}^{n+1} (q+1){i+1\choose
q+1}{n\choose i-1}= (q+1)\sum_i[{i-1\choose q+1}+2 {i-1\choose q}
+{i-1\choose q-1}]{n\choose i-1}$$ $$= (q+1)[{n\choose q+1}2^{n-q-1}
+ 2{n\choose q}2^{n-q}+{n\choose q-1}2^{n-q+1}].$$ Hence by (3.3) we
get
$$\beta_q(J_{n-1})= (q+1){n\choose q+1}2^{n-q-1} -2q{n\choose
q}2^{n-q}+(q-2){n\choose q-1}2^{n-q+1}.$$ On the other hand, by
Proposition 2.3,
$$\beta_q(I_{n-1}({\bf D})) = n {n-1\choose q}2^{n-q-1} +
(n-1){n\choose q-1}2^{n-q+1}.$$ Using $(3.1)$ this reduces to
$$\beta_q(I_{n-1}({\bf D})) = (q+1){n\choose q+1}2^{n-q-1} +
q{n+1\choose q}2^{n-q+1}-2{n\choose q-1}2^{n-q+1}.$$ By $(3.2)$
this is equal to $\beta_q(J_{n-1})$ computed above.\\
For $t=n-2,$ the proof is almost similar. More specifically, by
Lemma 3.3,
$$\beta_q(J_{n-2})=\sum_{\ell=q+1}^{6}{\ell-1\choose q}{n+ \ell
-4\choose n-3} + \sum_{j=q-5}^{n-3}{j+5\choose q}\sum_{\tau =1}^3
{n\choose j+\tau +2}{j+2\choose 3-\tau}.$$ Again we treat the case
$q \geq 6$ so that the first sum is zero. Using $i=j+3$ we get
$$\beta_q(J_{n-2})= \sum_{i=q-2}^{n}{i-1\choose 2}{i+2\choose q}{n\choose i} +
\sum_{i=q-2}^{n-1}(i-1){i+2\choose q}{n\choose i+1} + $$
$$\sum_{i=q-2}^{n-2}{i+2\choose q}{n\choose i+2}
=\sum_i[{i-1\choose 2}{i+2\choose q}+(i-2){i+1\choose q}+{i\choose
q}{n\choose i}.$$ Finally, with the same method as the case $t=n-1$
we arrive to
$$\beta_q(J_{n-2})= {q+2\choose 2}{n\choose q+2}2^{n-q-2} +4{q+1\choose 2}{n\choose
q+1}2^{n-q-1}+ $$ $$q(3q-4){n\choose q}2^{n-q}+2(q-1)(q-3){n\choose
q-1}2^{n-q+1}+{q-3\choose 2}{n\choose q-2}2^{n-q+2}.$$ By
Proposition 2.3 and manipulations as the previous case
$\beta_q(I_{n-2}({\bf D}))$ amounts to the same quantity.
$\hfill\square$
\paragraph{Acknowledgement.} The author thanks Dorin Pupesco for
suggesting him to compute the Betti numbers of the transversal
monomial ideals. He is grateful to Abdolhossein Hoorfar for his help
and comments on the combinatorics involved in the last theorem of
this paper.

\end{document}